\documentclass{article}
\usepackage{amssymb,amsmath,theorem,euscript}

\input{epsf}

\def\sm{\smallskip}

\newtheorem{theorem}{Theorem}

\newtheorem{lemma}[theorem]{Lemma}

\begin{document}

 \def\ov{\overline}
\def\wt{\widetilde}
\def\wh{\widehat}
 \newcommand{\rk}{\mathop {\mathrm {rk}}\nolimits}
\newcommand{\Aut}{\mathop {\mathrm {Aut}}\nolimits}
\newcommand{\Out}{\mathop {\mathrm {Out}}\nolimits}
\newcommand{\Abs}{\mathop {\mathrm {Abs}}\nolimits}
\renewcommand{\Re}{\mathop {\mathrm {Re}}\nolimits}
\renewcommand{\Im}{\mathop {\mathrm {Im}}\nolimits}
 \newcommand{\tr}{\mathop {\mathrm {tr}}\nolimits}
  \newcommand{\Hom}{\mathop {\mathrm {Hom}}\nolimits}
   \newcommand{\diag}{\mathop {\mathrm {diag}}\nolimits}
   \newcommand{\supp}{\mathop {\mathrm {supp}}\nolimits}
   \newcommand{\Ind}{\mathop {\mathrm {Ind}}\nolimits}
   \newcommand{\Rep}{\mathop {\mathrm {Rep}}\nolimits}
   \newcommand{\Irr}{\mathop {\mathrm {Irr}}\nolimits}
   \newcommand{\TR}{\mathop {\mathrm {TR}}\nolimits}
 \newcommand{\im}{\mathop {\mathrm {im}}\nolimits}
 \newcommand{\grad}{\mathop {\mathrm {grad}}\nolimits}
  \newcommand{\sgrad}{\mathop {\mathrm {sgrad}}\nolimits}
 \newcommand{\rot}{\mathop {\mathrm {rot}}\nolimits}
  \renewcommand{\div}{\mathop {\mathrm {div}}\nolimits}
  \newcommand{\sgn}{\mathop {\mathrm {sgn}}\nolimits}
   \newcommand{\dom}{\mathop {\mathrm {dom}}\nolimits}
  \newcommand{\sh}{\sinh}
  \newcommand{\ch}{\cosh}

\def\Br{\mathrm {Br}}
\def\Vir{\mathrm {Vir}}

 \def\Ham{\mathrm {Ham}}
\def\SL{\mathrm {SL}}
\def\Pol{\mathrm {Pol}}
\def\SU{\mathrm {SU}}
\def\GL{\mathrm {GL}}
\def\U{\mathrm U}
\def\OO{\mathrm O}
 \def\Sp{\mathrm {Sp}}
  \def\Ad{\mathrm {Ad}}
  \def\ad{\mathrm {ad}}
 \def\SO{\mathrm {SO}}
\def\SOS{\mathrm {SO}^*}
 \def\Diff{\mathrm{Diff}}
 \def\Vect{\mathfrak{Vect}}
\def\PGL{\mathrm {PGL}}
\def\PU{\mathrm {PU}}
\def\PSL{\mathrm {PSL}}
\def\Symp{\mathrm{Symp}}
\def\Cont{\mathrm{Cont}}
\def\End{\mathrm{End}}
\def\Mor{\mathrm{Mor}}
\def\Aut{\mathrm{Aut}}
 \def\PB{\mathrm{PB}}
\def\Fl{\mathrm {Fl}}
\def\Symm{\mathrm {Symm}} 
 \def\Herm{\mathrm {Herm}} 
  \def\SDiff{\mathrm {SDiff}} 
 
 \def\cA{\mathcal A}
\def\cB{\mathcal B}
\def\cC{\mathcal C}
\def\cD{\mathcal D}
\def\cE{\mathcal E}
\def\cF{\mathcal F}
\def\cG{\mathcal G}
\def\cH{\mathcal H}
\def\cJ{\mathcal J}
\def\cI{\mathcal I}
\def\cK{\mathcal K}
 \def\cL{\mathcal L}
\def\cM{\mathcal M}
\def\cN{\mathcal N}
 \def\cO{\mathcal O}
\def\cP{\mathcal P}
\def\cQ{\mathcal Q}
\def\cR{\mathcal R}
\def\cS{\mathcal S}
\def\cT{\mathcal T}
\def\cU{\mathcal U}
\def\cV{\mathcal V}
 \def\cW{\mathcal W}
\def\cX{\mathcal X}
 \def\cY{\mathcal Y}
 \def\cZ{\mathcal Z}
\def\0{{\ov 0}}
 \def\1{{\ov 1}}
 \def\frA{\mathfrak A}
 \def\frB{\mathfrak B}
\def\frC{\mathfrak C}
\def\frD{\mathfrak D}
\def\frE{\mathfrak E}
\def\frF{\mathfrak F}
\def\frG{\mathfrak G}
\def\frH{\mathfrak H}
\def\frI{\mathfrak I}
 \def\frJ{\mathfrak J}
 \def\frK{\mathfrak K}
 \def\frL{\mathfrak L}
\def\frM{\mathfrak M}
 \def\frN{\mathfrak N} \def\frO{\mathfrak O} \def\frP{\mathfrak P} \def\frQ{\mathfrak Q} \def\frR{\mathfrak R}
 \def\frS{\mathfrak S} \def\frT{\mathfrak T} \def\frU{\mathfrak U} \def\frV{\mathfrak V} \def\frW{\mathfrak W}
 \def\frX{\mathfrak X} \def\frY{\mathfrak Y} \def\frZ{\mathfrak Z} \def\fra{\mathfrak a} \def\frb{\mathfrak b}
 \def\frc{\mathfrak c} \def\frd{\mathfrak d} \def\fre{\mathfrak e} \def\frf{\mathfrak f} \def\frg{\mathfrak g}
 \def\frh{\mathfrak h} \def\fri{\mathfrak i} \def\frj{\mathfrak j} \def\frk{\mathfrak k} \def\frl{\mathfrak l}
 \def\frm{\mathfrak m} \def\frn{\mathfrak n} \def\fro{\mathfrak o} \def\frp{\mathfrak p} \def\frq{\mathfrak q}
 \def\frr{\mathfrak r} \def\frs{\mathfrak s} \def\frt{\mathfrak t} \def\fru{\mathfrak u} \def\frv{\mathfrak v}
 \def\frw{\mathfrak w} \def\frx{\mathfrak x} \def\fry{\mathfrak y} \def\frz{\mathfrak z} \def\frsp{\mathfrak{sp}}
 \def\bfa{\mathbf a} \def\bfb{\mathbf b} \def\bfc{\mathbf c} \def\bfd{\mathbf d} \def\bfe{\mathbf e} \def\bff{\mathbf f}
 \def\bfg{\mathbf g} \def\bfh{\mathbf h} \def\bfi{\mathbf i} \def\bfj{\mathbf j} \def\bfk{\mathbf k} \def\bfl{\mathbf l}
 \def\bfm{\mathbf m} \def\bfn{\mathbf n} \def\bfo{\mathbf o} \def\bfp{\mathbf p} \def\bfq{\mathbf q} \def\bfr{\mathbf r}
 \def\bfs{\mathbf s} \def\bft{\mathbf t} \def\bfu{\mathbf u} \def\bfv{\mathbf v} \def\bfw{\mathbf w} \def\bfx{\mathbf x}
 \def\bfy{\mathbf y} \def\bfz{\mathbf z} \def\bfA{\mathbf A} \def\bfB{\mathbf B} \def\bfC{\mathbf C} \def\bfD{\mathbf D}
 \def\bfE{\mathbf E} \def\bfF{\mathbf F} \def\bfG{\mathbf G} \def\bfH{\mathbf H} \def\bfI{\mathbf I} \def\bfJ{\mathbf J}
 \def\bfK{\mathbf K} \def\bfL{\mathbf L} \def\bfM{\mathbf M} \def\bfN{\mathbf N} \def\bfO{\mathbf O} \def\bfP{\mathbf P}
 \def\bfQ{\mathbf Q} \def\bfR{\mathbf R} \def\bfS{\mathbf S} \def\bfT{\mathbf T} \def\bfU{\mathbf U} \def\bfV{\mathbf V}
 \def\bfW{\mathbf W} \def\bfX{\mathbf X} \def\bfY{\mathbf Y} \def\bfZ{\mathbf Z} \def\bfw{\mathbf w}
 \def\R {{\mathbb R }} \def\C {{\mathbb C }} \def\Z{{\mathbb Z}} \def\H{{\mathbb H}} \def\K{{\mathbb K}}
 \def\N{{\mathbb N}} \def\Q{{\mathbb Q}} \def\A{{\mathbb A}} \def\T{\mathbb T} \def\P{\mathbb P} \def\G{\mathbb G}
 \def\bbA{\mathbb A} \def\bbB{\mathbb B} \def\bbD{\mathbb D} \def\bbE{\mathbb E} \def\bbF{\mathbb F} \def\bbG{\mathbb G}
 \def\bbI{\mathbb I} \def\bbJ{\mathbb J} \def\bbL{\mathbb L} \def\bbM{\mathbb M} \def\bbN{\mathbb N} \def\bbO{\mathbb O}
 \def\bbP{\mathbb P} \def\bbQ{\mathbb Q} \def\bbS{\mathbb S} \def\bbT{\mathbb T} \def\bbU{\mathbb U} \def\bbV{\mathbb V}
 \def\bbW{\mathbb W} \def\bbX{\mathbb X} \def\bbY{\mathbb Y} \def\kappa{\varkappa} \def\epsilon{\varepsilon}
 \def\phi{\varphi} \def\le{\leqslant} \def\ge{\geqslant}

\def\UU{\bbU}
\def\Mat{\mathrm{Mat}}
\def\tto{\rightrightarrows}

\def\F{\mathbf{F}}

\def\Gms{\mathrm {Gms}}
\def\Ams{\mathrm {Ams}}
\def\Isom{\mathrm {Isom}}

\def\Gr{\mathrm{Gr}}

\def\graph{\mathrm{graph}}

\def\O{\mathrm{O}}

\def\la{\langle}
\def\ra{\rangle}


 \def\ov{\overline}
\def\wt{\widetilde}

\renewcommand{\Re}{\mathop {\mathrm {Re}}\nolimits}
\def\Br{\mathrm {Br}}

 \def\Isom{\mathrm {Isom}}
 \def\Free{\mathrm {Free}}
 \def\Hier{\mathrm {Hier}}
  \def\Comm{\mathrm {Comm}}
\def\SL{\mathrm {SL}}
\def\SU{\mathrm {SU}}
\def\GL{\mathrm {GL}}
\def\U{\mathrm U}
\def\OO{\mathrm O}
 \def\Sp{\mathrm {Sp}}
  \def\GLO{\mathrm {GLO}}
 \def\SO{\mathrm {SO}}
\def\SOS{\mathrm {SO}^*}
 \def\Diff{\mathrm{Diff}}
 \def\Vect{\mathfrak{Vect}}
\def\PGL{\mathrm {PGL}}
\def\PU{\mathrm {PU}}
\def\PSL{\mathrm {PSL}}
\def\Symp{\mathrm{Symp}}
\def\ASymm{\mathrm{Asymm}}
\def\Asymm{\mathrm{Asymm}}
\def\Gal{\mathrm{Gal}}
\def\End{\mathrm{End}}
\def\Mor{\mathrm{Mor}}
\def\Aut{\mathrm{Aut}}
 \def\PB{\mathrm{PB}}
 \def\cA{\mathcal A}
\def\cB{\mathcal B}
\def\cC{\mathcal C}
\def\cD{\mathcal D}
\def\cE{\mathcal E}
\def\cF{\mathcal F}
\def\cG{\mathcal G}
\def\cH{\mathcal H}
\def\cJ{\mathcal J}
\def\cI{\mathcal I}
\def\cK{\mathcal K}
 \def\cL{\mathcal L}
\def\cM{\mathcal M}
\def\cN{\mathcal N}
 \def\cO{\mathcal O}
\def\cP{\mathcal P}
\def\cQ{\mathcal Q}
\def\cR{\mathcal R}
\def\cS{\mathcal S}
\def\cT{\mathcal T}
\def\cU{\mathcal U}
\def\cV{\mathcal V}
 \def\cW{\mathcal W}
\def\cX{\mathcal X}
 \def\cY{\mathcal Y}
 \def\cZ{\mathcal Z}
\def\0{{\ov 0}}
 \def\1{{\ov 1}}
 
 \def\frA{\mathfrak A}
 \def\frB{\mathfrak B}
\def\frC{\mathfrak C}
\def\frD{\mathfrak D}
\def\frE{\mathfrak E}
\def\frF{\mathfrak F}
\def\frG{\mathfrak G}
\def\frH{\mathfrak H}
\def\frI{\mathfrak I}
 \def\frJ{\mathfrak J}
 \def\frK{\mathfrak K}
 \def\frL{\mathfrak L}
\def\frM{\mathfrak M}
 \def\frN{\mathfrak N} \def\frO{\mathfrak O} \def\frP{\mathfrak P} \def\frQ{\mathfrak Q} \def\frR{\mathfrak R}
 \def\frS{\mathfrak S} \def\frT{\mathfrak T} \def\frU{\mathfrak U} \def\frV{\mathfrak V} \def\frW{\mathfrak W}
 \def\frX{\mathfrak X} \def\frY{\mathfrak Y} \def\frZ{\mathfrak Z} \def\fra{\mathfrak a} \def\frb{\mathfrak b}
 \def\frc{\mathfrak c} \def\frd{\mathfrak d} \def\fre{\mathfrak e} \def\frf{\mathfrak f} \def\frg{\mathfrak g}
 \def\frh{\mathfrak h} \def\fri{\mathfrak i} \def\frj{\mathfrak j} \def\frk{\mathfrak k} \def\frl{\mathfrak l}
 \def\frm{\mathfrak m} \def\frn{\mathfrak n} \def\fro{\mathfrak o} \def\frp{\mathfrak p} \def\frq{\mathfrak q}
 \def\frr{\mathfrak r} \def\frs{\mathfrak s} \def\frt{\mathfrak t} \def\fru{\mathfrak u} \def\frv{\mathfrak v}
 \def\frw{\mathfrak w} \def\frx{\mathfrak x} \def\fry{\mathfrak y} \def\frz{\mathfrak z} \def\frsp{\mathfrak{sp}}
 \def\bfa{\mathbf a} \def\bfb{\mathbf b} \def\bfc{\mathbf c} \def\bfd{\mathbf d} \def\bfe{\mathbf e} \def\bff{\mathbf f}
 \def\bfg{\mathbf g} \def\bfh{\mathbf h} \def\bfi{\mathbf i} \def\bfj{\mathbf j} \def\bfk{\mathbf k} \def\bfl{\mathbf l}
 \def\bfm{\mathbf m} \def\bfn{\mathbf n} \def\bfo{\mathbf o} \def\bfp{\mathbf p} \def\bfq{\mathbf q} \def\bfr{\mathbf r}
 \def\bfs{\mathbf s} \def\bft{\mathbf t} \def\bfu{\mathbf u} \def\bfv{\mathbf v} \def\bfw{\mathbf w} \def\bfx{\mathbf x}
 \def\bfy{\mathbf y} \def\bfz{\mathbf z} \def\bfA{\mathbf A} \def\bfB{\mathbf B} \def\bfC{\mathbf C} \def\bfD{\mathbf D}
 \def\bfE{\mathbf E} \def\bfF{\mathbf F} \def\bfG{\mathbf G} \def\bfH{\mathbf H} \def\bfI{\mathbf I} \def\bfJ{\mathbf J}
 \def\bfK{\mathbf K} \def\bfL{\mathbf L} \def\bfM{\mathbf M} \def\bfN{\mathbf N} \def\bfO{\mathbf O} \def\bfP{\mathbf P}
 \def\bfQ{\mathbf Q} \def\bfR{\mathbf R} \def\bfS{\mathbf S} \def\bfT{\mathbf T} \def\bfU{\mathbf U} \def\bfV{\mathbf V}
 \def\bfW{\mathbf W} \def\bfX{\mathbf X} \def\bfY{\mathbf Y} \def\bfZ{\mathbf Z} \def\bfw{\mathbf w}

 \def\R {{\mathbb R }} \def\C {{\mathbb C }} \def\Z{{\mathbb Z}} \def\H{{\mathbb H}} \def\K{{\mathbb K}}
 \def\N{{\mathbb N}} \def\Q{{\mathbb Q}} \def\A{{\mathbb A}} \def\T{\mathbb T} \def\P{\mathbb P} \def\G{\mathbb G}
 \def\bbA{\mathbb A} \def\bbB{\mathbb B} \def\bbD{\mathbb D} \def\bbE{\mathbb E} \def\bbF{\mathbb F} \def\bbG{\mathbb G}
 \def\bbI{\mathbb I} \def\bbJ{\mathbb J} \def\bbL{\mathbb L} \def\bbM{\mathbb M} \def\bbN{\mathbb N} \def\bbO{\mathbb O}
 \def\bbP{\mathbb P} \def\bbQ{\mathbb Q} \def\bbS{\mathbb S} \def\bbT{\mathbb T} \def\bbU{\mathbb U} \def\bbV{\mathbb V}
 \def\bbW{\mathbb W} \def\bbX{\mathbb X} \def\bbY{\mathbb Y} \def\kappa{\varkappa} \def\epsilon{\varepsilon}
 \def\phi{\varphi} \def\le{\leqslant} \def\ge{\geqslant}

\def\UU{\bbU}
\def\Mat{\mathrm{Mat}}
\def\tto{\rightrightarrows}

\def\Gr{\mathrm{Gr}}

\def\B{\bfB} 

\def\graph{\mathrm{graph}}

\def\O{\mathrm{O}}

\def\la{\langle}
\def\ra{\rangle}

	\begin{center}
		\Large\bf
		Fourier transform on the Lobachevsky plane and operational calculus

		\bigskip
		
		\large\sc
		Yu.~A. Neretin%
		\footnote{The work is supported by the grant FWF, P31591.}
	\end{center}

{\small The classical Fourier transform
	on the line sends the operator of multiplication by   $x$
to $i\frac{d}{d\xi}$ and the operator
of differentiation  $\frac{d}{d x}$ to the multiplication by  $-i\xi$.
For the Fourier transform on the Lobachevsky plane we establish
a similar correspondence for a certain family of differential operators. 
It appears that differential operators on the Lobachevsky plane correspond
to differential-difference operators in the Fourier-image, where
shift operators act in the imaginary direction, i.e., a direction transversal to the integration
contour in the Plancherel formula. }

\bigskip

{\bf 1. Lobachevsky plane.} Consider the complex plane
with the coordinate
 $z=x+iy$ and the upper half-plane
 $\Lambda$ consisting of points with $y=\Im z>0$.
Consider the group $\SL(2,\R)$ of real matrices  $g=\begin{pmatrix}
a&b\\c&d
\end{pmatrix}$ with determinant 1. It acts on  $\Lambda$
by linear-fractional maps
$$
z\mapsto z^{[g]}:=\frac{b+xd}{a+xc}.
$$
Then $\Lambda$ becomes a homogeneous space ({\it the Lobachevsky plane})
$$
\Lambda=\SO(2)\setminus \SL(2,\R),
$$
the subgroup
 $\SO(2)$ consisting of matrices  $\begin{pmatrix}
\cos \phi&\sin\phi\\-\sin\phi&\cos\phi
\end{pmatrix}$ is the stabilizer of the point  $i\in \Lambda$.
Notice that the transformation corresponding to the matrix
 $\begin{pmatrix}
-1&0\\0&-1\end{pmatrix}$ is identical, so in fact we have an action
of the quotient group
$$
\PSL(2,\R):=\SL(2,\R)/\{\pm 1 \}.
$$
The transformations
 $z\to z^{[g]}$ induce transformation of a space of functions
on
 $\Lambda$,
$$
R(g)f(z):=f(z^{[g]}).
$$

It is easy to verify that the measure
$$
d\mu(x,y)=d\mu(z)=\frac{dx\,dy}{4\,y^2}
$$
is invariant with respect to the transformations
$z\mapsto z^{[g]}$. Therefore the transformations  $R(g)$
determine a unitary representation of the group
 $\SL(2,\R)$ in the space  $L^2$ on $\Lambda$ with respect to the measure $\mu$.

\sm

{\bf 2. The principal series of unitary representations of $\PSL(2,\R)$.}
See \cite{GGV}, Chapter VII.
For $\tau\in \C$ we define a representation   $T_\tau$ 
of the group
$\SL(2,\R)$ in a space of functions on the line by the formula
$$
T_\tau\begin{pmatrix}
a&b\\c&d
\end{pmatrix} f(x)
=f\Bigl(\frac{b+xd}{a+xc}\Bigr) |a+cx|^{2\tau}.
 $$
If $\Re\tau=-\frac 12$, then this representation
is unitary in  $L^2(\R)$. Moreover, representations 
$T_{-1/2+is}$ and $T_{-1/2-is}$ are equivalent. Representations $T_\tau$ with $\Re\tau=-1/2$
are called {\it representations of the unitary principal series}.

If $\Re\tau\ne-1/2$, then we need some care
to define a space of representation. 
For us it will be convenient the following version.
Denote by
$C^\infty_\tau(\R)$ the space of $C^\infty$-smooth functions  $f$ on $\R$,
satisfying the following additional condition%
\footnote{At a first glance this condition seems awkward.
	In fact it is more natural to consider the space of smooth functions
 (or spaces of smooth sections of linear bundles) on the circle 
(the projective line) $\R\cup \infty$. Passing to a space
of functions on the line we cut the circle, for this reason
we must impose  conditions of gluing at infinity.}:
 A function $f(1/x)|x|^{2\tau}$ has a removable singularity at 0 and becomes 
 $C^\infty$-smooth after the removing.
This condition provides invariance of the space 
 $C^\infty_\tau(\R)$ with respect to the transformations  $T_\tau(g)$. .

Representations obtained in this way are called
 {\it representations of the principal  {\rm(}nonunitary{\rm)} series.}
 Notice that for
 $\tau\notin  \Z$ they are irreducible,  $T_{-\lambda+1/2}$ is equivalent to $T_{\lambda+1/2}$.

\sm

{\bf 3. The Fourier transform on the Lobachevsky plane.}
Denote by
 $C_c(\Lambda)$ the space of smooth compactly supported functions
  on $\Lambda$.
Set
$$
K(\tau;z,x):=\left(\frac{2i(x-z)(x-\ov z)}{z-\ov z}  \right)^\tau,\qquad\text{where $z\in\Lambda$, $x\in\R$.}
$$	
For $f\in C_c(\Lambda)$ we assign the function
  (the {\it Fourier transform})
$Jf$ on $\C\times \R$ by the formula 
\begin{equation}
Jf(\tau; x)=\int_{\Lambda} K(\tau;z,x) f(z) \,d \mu(z).
\end{equation}	
It is easy to verify that
$$
JR(g)f(\tau;x)=f(\tau;x^{[g]})|a+xc|^{2\tau} = T_\tau(g)f(\tau;x).
$$
This can be easily derived from the formula
$$
u^{[g]}-v^{[g]}=\frac{u-v}{(a+uc)(a+vc)},\qquad \text{where $g=\begin{pmatrix}
	a&b\\c&d
	\end{pmatrix}$, $\det g=1$.}
$$ 

The image $\cP(\C \times \R)$ of the space  $C_c(\Lambda)$
admits a precise description (a '{\it Paley--Wiener theorem}',
see S.~Helgason \cite{Hel}, \S I.4, see also \cite{Ter}, \S 3.2).
 A function $\Phi(\tau;x)$ is contained in $\cP(\C \times \R)$
   if it satisfies the following conditions
 0) -- 2):

\sm 

0) The functions $\Phi(\tau;x)$ and $\Phi(\tau;1/x)|x|^{2\tau}$ are $C^\infty$-smooth 
on $\C\times \R$ and are holomorphic in  $\tau$.

\sm 

1) Fix $A>1$ (a choice of $A$ has no matter). For a function $\Phi$
there exists
$R$ such that for any   $N$ the following condition holds
\begin{align*}
\sup_{\tau\in \C, x\in [-A,A]}
|\Phi(\tau;x)|\, e^{-R|\Re \tau|} (1+|\tau|)^N< \infty;
\\
\sup_{\tau\in \C, x\in [-A,A]}
\bigl|\Phi(\tau;1/x)|x|^{2\tau}\bigr|\, e^{-R|\Re \tau|} (1+|\tau|)^N< \infty.
\end{align*}

2) The following identity holds 
\begin{multline*}
\int_{-\infty}^\infty
\Phi(-1/2+\lambda;x) K(-1/2-\lambda;z,x)\,dx=\\=
\int_{-\infty}^\infty
\Phi(-1/2-\lambda;x) K(-1/2+\lambda;z,x)\,dx.
\end{multline*}

Notice that condition
 1) is similar to the condition of the usual 'Paley--Wiener theorem'  in the sense
 of L.~Schwartz. Condition 2 means a kind of 'evenness',
the functions $\phi^\lambda(x):=\Phi(-1/2+\lambda;x)$ and $\phi^{-\lambda}(x):=\Phi(-1/2-\lambda;x)$
determine one another. They are not equal, but are related by an integral condition.

Another fundamental statement about the Fourier transform is the  
{\it Plancherel theorem}.
For $f_1$, $f_2\in C_c^\infty(\Lambda)$ the following identity holds
\begin{multline}
\int_{\Lambda} f_1(z) \ov{f_2(z)}\,d\mu(z)=\\=
\int_0^\infty \int_{-\infty}^\infty Jf_1(-1/2+is;x)\ov {Jf_2(-1/2+is)}\,\, s\, \frac{\sh (\pi s)}{\ch (\pi s)}\,
\,dx\, ds.
\label{eq:sh}
\end{multline}
Moreover, the operator
 $J$ extends to a unitary operator from $L^2(\Lambda,d\mu)$ to
 $L^2$ on the domain  $s\ge 0$, $x\in \R$ with respect to the measure 
$ s \frac{\sh (\pi s)}{\ch (\pi s)}\,dx\, ds$. 

It is quite easy to write an explicit form of the intertwining operator
 $J$ providing the  spectral decomposition of
 the representation
 $R(g)$. However an explicit expression
 $s\frac{\sh (\pi s)}{\ch (\pi s)}\, ds$ for the spectral measure is a relatively
 delicate fact.
 This statement was formulated by Ferdinand Mehler
 \cite{Meh} in 1881 without proof. Various proofs were published
 by H.~Weyl  \cite{Wey} in 1910 (as a very particular case
 of spectral theory of differential operators), V.~A.~Fock  \cite{Fock} in 1943,
M.~N.~Olevsky \cite{Ole} in 1949. See, also \cite{Koo}, \cite{Hel}, \S.I.4, \cite{Ter}, \S 3.2.

\sm 

{\bf 4. The correspondence of differential operators.}
Our 'Fourier transform' is one of simplest
representatives of numerous 'Fourier transforms' or 'Plancherel decompositions'
in noncommutative harmonic analysis. 
Sometimes such transformations are relatively strong tool
for theory of special functions, more often they are more points on
the boundary of a 
knowable world than standpoints of a new life.
In particular there are not too much functions, for which 
Fourier transform admits an explicit evaluation.

The author in \cite{Ner-over} observed that for
one transformation of such type there is a kind of an operational
calculus. Later there appeared several
works of V.~F.~Molchanov \cite{Mol-faa1}--\cite{Mol-faa2} 
and the author \cite{Ner-gln}--\cite{Ner-gl}  on this subject.
See a wider discussion in \cite{Ner-sl2}, \cite{Ner-after}.

We start from an obvious statement.
Since the Fourier transform commutes with the action of the Lie group
 $\SL(2,\R)$, it commutes with the action of the Lie algebra
 $\frs\frl(2)$. This implies the following correspondence
 of differential operators and their images under the Fourier transform:
\begin{align}
\frac{\partial}{\partial z}+\frac{\partial}{\partial \ov z}\quad
&\longleftrightarrow\quad  \frac{\partial}{\partial x};
\label{eq:corr01}
\\
z\frac{\partial}{\partial z}+ \ov z \frac{\partial}{\partial \ov z}\quad
&\longleftrightarrow\quad x \frac{\partial}{\partial x}-\tau;
\label{eq:corr02}
\\
z^2\frac{\partial}{\partial z}+ \ov z^2 \frac{\partial}{\partial \ov z}\quad
&\longleftrightarrow\quad x^2 \frac{\partial}{\partial x}-2\tau x;
\label{eq:corr03}
\end{align}
These formulas can be easily verified in a straightforward way.

Let $\Phi(\tau,x)$ be a function on $\C\times \R$. We define operators 
\begin{equation*}
T_+ \Phi(\tau; x):= \Phi(\tau+1; x),\qquad T_- \Phi(\tau; x):= \Phi(\tau-1; x)
\end{equation*}

\begin{theorem}
	\label{th:1}
	The Fourier transform establishes  the following correspondence between differential operators
	in the space
 $C_c^\infty(\Lambda)$
 and differential-difference operators in the space
 $\cP(\C \times \R)$:
 	\begin{align}
 \frac {1}{z-\ov z}
 \qquad& \longleftrightarrow \qquad  \frac 1{4i (1+\tau)(1+2\tau)}
 \, \frac{\partial^2}{\partial x^2} T_+ -
 \frac{2i\tau}{2(1+2\tau)}\, T_-\,;
 \label{eq:corr10}
 \\	
 \frac{\partial}{\partial z}- \frac{\partial}{\partial \ov z}
 \qquad &\longleftrightarrow \qquad 
 \frac{2+\tau}{2i (1+\tau)(1+2\tau)}\,\frac{\partial^2 }{\partial x^2}\, T_+ +
 \frac{2i\tau(-1+\tau)}{1+2\tau}\, T_-\,;
 \label{eq:corr11}
 \\
 z	\frac{\partial}{\partial z}- \ov z\frac{\partial}{\partial \ov z}
 \qquad &\longleftrightarrow \qquad
 \frac{(2+\tau)}{2i(1+\tau)(1+2\tau)}\, x \frac{\partial^2 }{\partial x^2}\, T_+-
 \frac{(2+\tau)}{2i(1+\tau)} \,\frac{\partial }{\partial x}\, T_+
 \nonumber
 +\\&
 \qquad\qquad\qquad\qquad\qquad
 +\frac{2i(-1+\tau)\tau}{(1+2\tau)}\,x T_-\,;
 \label{eq:corr12}
 \\
 z^2	\frac{\partial}{\partial z}- \ov z^{\,2}\frac{\partial}{\partial \ov z}
 \qquad &\longleftrightarrow \qquad
 \frac{(2+\tau)}{2i(1+\tau)(1+2\tau)}\, x^2 \frac{\partial^2 }{\partial x^2}\,T_+
 -
 \frac{2(2+\tau)}{2i(1+\tau)} \, x \frac{\partial }{\partial x}\,T_+ +
 \nonumber
 \\&
 \qquad\qquad\qquad
 +\frac{2(2+\tau)}{2i} \,T_+
 + \frac{2i\tau(-1+\tau)}{(1+2\tau)}\,x^2 T_-\,.
 \label{eq:corr13}
 \end{align}
\end{theorem}

{\sc Remarks.}
1) The vector fields $\frac{\partial}{\partial z}$, $z\frac{\partial}{\partial z}$, 
$z^2\frac{\partial}{\partial z}$ form a Lie algebra isomorphic to $\frs\frl(2)$.
The same holds for vector fields  $\frac{\partial}{\partial \ov z}$, $\ov z\frac{\partial}{\partial\ov z}$, 
$\ov z^2\frac{\partial}{\partial\ov z}$ and these Lie algebras of vector fields commute.
Our formulas provides us the Fourier-images for all elements of $\frs\frl(2)\oplus \frs\frl(2)$. 

\sm

2) The shift operators in these formulas have imaginary
direction with respect to the surface of integration in
(\ref{eq:sh}). This phenomenon take place for all known problems of this
kind related to semisimple Lie groups  (see \cite{Ner-over},
\cite{Mol-faa1}--\cite{Mol-faa2}, \cite{Ner-gln}--\cite{Ner-gl}).
Appearance of differential operators of high order in Fourier-images of simplest operators also
is usual, see  \cite{Mol-faa1}, \cite{Ner-gln} (our group is 'small'
therefore the order of operators is only 2). Although shifts in imaginary
direction  $\Phi(t)\mapsto \Phi(t+i)$ seem strange from the point of view
of $L^2$-theory, related problems of spectral theory are quite reasonable
 (see \cite{Gro}, \cite{Ner-difference},
\cite{MN}). Notice that there are problems of this kind 
for collections of commuting operators (see I.~Cherednik
 \cite{Che}, see also J. van Diejen, E.~Emez \cite{DE}).
\hfill $\lozenge$

\sm 

{\bf 5. The proof of Theorem \ref{th:1}.}
The first correspondence (\ref{eq:corr10}) is equivalent to the identity
\begin{multline}
\frac{1}{z-\ov z} K(\tau;z,x)=\Bigl( - \frac 1{4i (1+\tau)(1+2\tau)}
\, \frac{\partial^2}{\partial x^2} T_+ 
+\\+
\frac{2i\tau}{2(1+2\tau)}\, T_-\Bigr) K(\tau;z,x).
\label{eq:id1} 
\end{multline}
Dividing both sides by
 $K(\tau;z,x)$, we observe that it is sufficient to verify
 the following identity for rational functions:
\begin{multline}
\frac{1}{z-\ov z} =K(\tau;z,x)^{-1}\Bigl( - \frac 1{4i (1+\tau)(1+2\tau)}
\, \frac{\partial^2}{\partial x^2} T_+ 
+\\+
\frac{2i\tau}{2(1+2\tau)}\, T_-\Bigr) K(\tau;z,x),
\label{eq:id2} 
\end{multline}

For evaluation of the Fourier-image of the second operator
 $\partial/\partial z-
\partial/\partial \ov z$ we write
$$
\Bigl[J\circ \Bigl(\frac{\partial}{\partial z}-
\frac{\partial}{\partial \ov z}\Bigr)\Bigr] f(x)=
\int_{\Lambda} K(\tau;z,x) \Bigl(\frac{\partial}{\partial z}-
\frac{\partial}{\partial \ov z}\Bigr) f(z)\frac{dz\,d\ov z}{2i(z-\ov z)^2}.
$$
Integrating by parts we come to
$$
\int_{\Lambda} \Bigl(-\frac{\partial}{\partial z}+
\frac{\partial}{\partial \ov z}+\frac4{z-\ov z}\Bigr) K(\tau;z,x) \cdot f(z)\frac{dz\,d\ov z}{2i(z-\ov z)^2}.
$$
We must get the operator from the right-hand side of
 (\ref{eq:corr11}), i.~e., we must verify the identity 
\begin{multline}
\Bigl(-\frac{\partial}{\partial z}+\frac{\partial}{\partial \ov z}+\frac4{z-\ov z}\Bigr) K(\tau;z,x)
=\\=
\Bigl(\frac{2+\tau}{2i (1+\tau)(1+2\tau)}\,\frac{\partial^2 }{\partial x^2}\, T_+ +
\frac{2i\tau(-1+\tau)}{1+2\tau}\, T_-\Bigr)K(\tau;z,x).
\label{eq:second}
\end{multline}
Dividing both sides by
 $K(\tau;z,x)$, we as in (\ref{eq:id2}), 
 get rational functions in both sides of the equality.
 Such an identity can be easily verified with
 Mathematica or Maple.

Let us explain how to verify formula
 (\ref{eq:id2}) by hands.  Denote
$$
A:= \frac{2i\,T_- K(\tau;z,x)}
{K(\tau;z,x)}, \qquad B:= \frac{- \frac 1{2i}
\, \frac{\partial^2}{\partial x^2} T_+K(\tau;z,x)}{K(\tau;z,x)}.
$$
Then
\begin{equation*}
A=\frac{z-\ov z}{(x-z)(x-\ov z)}=\frac{1}{x-z}-\frac{1}{x-\ov z}\,,
\end{equation*}
and
\begin{multline*}
\frac 1{\tau+1}B=
\Bigl(\frac{\tau}{(x-z)^2}+ \frac{\tau}{(x-\ov z)^2}+\frac{2(\tau+1)}{(x-z)(x-\ov z)}\Bigr)
\frac{(x-z)(x-\ov z)}{z-\ov z}=\\
=\frac{\tau}{x-z}- \frac{\tau}{x-\ov z}-\frac{2(1+\tau)}{z-\ov z}.
\end{multline*}

For a verification of (\ref{eq:id2}) it is sufficient to evaluate
coefficients at prime fractions. 

\sm

The identity  (\ref{eq:corr11}) can be verified in a similar way.
Additionally, we write
\begin{equation*}
\frac{\Bigl(\frac{\partial}{\partial z}- \frac{\partial}{\partial \ov z}\Bigr)K(\tau;z,x)}
{K(\tau;z,x)}=-\frac{\tau}{x-z}+\frac{\tau}{x-\ov z}-\frac{2\tau}{z-\ov z},
\end{equation*}
and verify coefficients at prime fractions.

Now we can take the commutator of the
 (\ref{eq:corr11}) with the operator 
(\ref{eq:corr03}), this gives the correspondence  (\ref{eq:corr12}).
Evaluating the commutator of 
(\ref{eq:corr12}) with (\ref{eq:corr03}), we come to the correspondence (\ref{eq:corr13}).
\hfill $\square$

\sm

{\bf 6. An algebra of operators whose Fourier-images admit evaluations.}
\begin{figure}
$$
\epsfbox{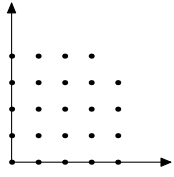}
$$
\caption{A set $\Omega_n$.}
\label{fig:1}
\end{figure}
By $\Omega_n$ we denote the set of all integer points
$(p,q)$ such that
\begin{equation}
0\le p\le n,\quad 0\le q\le n, \quad (p,q)\ne(n,n),
\end{equation}
see Fig. \ref{fig:1}.
Consider the space
 $\C(z,\ov z)$ of rational functions of variables $z$ and $\ov z$.
Denote by $\cA_n\subset \C(z,\ov z)$ the subspace with the basis 
\begin{equation}
\frac{z^p\,\ov z^{\, q}}{(z-\ov z)^n},\quad \text{where $(p,q)\in\Omega_n$.}
\label{eq:basis}
\end{equation} 
Denote by
 $\cA\subset \C(z,\ov z)$ the sum of all subspaces   $\cA_n$.
 Clearly $\cA_n$ is a subalgebra in $\C(z,\ov z)$.
 Consider the space of all differential operators on Lobachevsky plane,
 denote by
  $\cB$ the subspace consisting of all 
  operators of the form
\begin{multline}
\label{eq:cB}
\sum\limits_{\alpha,\beta,\gamma,\alpha',\beta',\gamma'\ge 0}
V_{\alpha,\beta,\gamma,\alpha',\beta',\gamma'}(z,\ov z) \,
\Bigl(z^2\frac{\partial}{\partial z}\Bigr)^{\alpha} \Bigl(z\frac{\partial}{\partial z}\Bigr)^{\beta}
\Bigl(\frac{\partial}{\partial z}\Bigr)^{\gamma}
\times \\ \times
\Bigl(\ov z^2\frac{\partial}{\partial \ov z}\Bigr)^{\alpha'} 
\Bigl(\ov z\frac{\partial}{\partial \ov z}\Bigr)^{\beta'}
\Bigl(\frac{\partial}{\partial \ov z}\Bigr)^{\gamma'},
\end{multline}
where $V_{\alpha,\beta,\gamma,\alpha',\beta',\gamma'}(z,\ov z)\in \cA$.

\begin{theorem}
\label{th:2}
The space
 $\cB$ is a subalgebra in the algebra of all differential operators.
 This subalgebra is generated by the operators
\begin{equation}
\frac{1}{z-\ov z},\,\,\,z^2\frac{\partial}{\partial z},\,\,\, z\frac{\partial}{\partial z},\,\,\,
\frac{\partial}{\partial z},\,\,\, \ov z^2\frac{\partial}{\partial\ov z},\,\,\,
\ov z\frac{\partial}{\partial\ov z},\,\,\, \frac{\partial}{\partial\ov z}.
\label{eq:generators}
\end{equation}
\end{theorem}

We also can say that $\cB$ is the space of all expressions of the form

\begin{theorem}
\label{th:3}
For any element of the algebra
 $\cB$ the corresponding operator in the Fourier-image
 has the form
\begin{equation}
\sum\limits_{p\ge 0,q\ge 0,r\in \Z}
U_{p,q,r}(\tau) \,x^p\frac{\partial^q}{\partial x^q} T_+^r,
\label{eq:sumU}
\end{equation}
where $U_{p,q,r}(\tau)$ are rational functions in  $\tau$
with possible poles at points $\tau\in\Z/2$.
\end{theorem}

Theorem \ref{th:3} is a straightforward corollary of Theorem
\ref{th:2} and statements, which were established above. 
Indeed, the space of operators of the form  (\ref{eq:sumU})
is an algebra, denote it by
 $\cD$. In virtue of formulas  (\ref{eq:corr01})--(\ref{eq:corr13}),
 the Fourier-images of all generators
 (\ref{eq:generators}) of the algebra  $\cB$ are contained in $\cD$.

\sm 

{\bf 7. The proof of Theorem \ref{th:2}.}

\begin{lemma}
\label{l:1}
{\rm a)} The space $\cA$ is invariant with respect to the operators 
\begin{equation}
\frac\partial{\partial z},\,\, z\frac\partial{\partial z},\,\,z^2\frac\partial{\partial z},\,\,
\frac\partial{\partial \ov z},\,\,\ov z\frac\partial{\partial\ov z},
\,\,\ov z^2\frac\partial{\partial \ov z}.
\label{eq:6-operatorov}
\end{equation}

{\rm b)} The function $1/(z-\ov z)$ is cyclic in the space   $\cA$
with respect to this family of operators.
In other words, apply all possible products of such operators to
 $1/(z-\ov z)$ and consider the subspace
 $\wt \cA\subset \cA$ spanned by such functions. Then
$\wt \cA=\cA$.  
\end{lemma}

{\sc Proof.}
Denote by $\cA_{[n]}\subset \cA$ the sum of all subspaces  $\cA_j$ with $j\le n$. 

We prove the statement by induction. Notice that
$$
\Bigl(z^2\frac\partial{\partial z}+\ov z^2\frac\partial{\partial \ov z}\Bigr) 
\frac{1}{z-\ov z}=\frac{z+\ov z}{z-\ov z}\in\wt \cA.
$$
Since
$$
1=\frac{z-\ov z}{z-\ov z},
$$
we get  
$$
\frac{z}{z-\ov z},\,\, \frac{\ov z}{z-\ov z}\in\wt\cA.
$$
Also, notice that applying any operator
 (\ref{eq:6-operatorov}) to $1/(z-\ov z)$, we get an element of the space
 $\cA_2$.

Fix
 $n$. Assume that  $\cA_{[n]}\subset \wt \cA$.
 Let apply 6 operators (\ref{eq:6-operatorov}) to all functions 
$z^p\,\ov z^q/(x-y)^n\in\cA_n$ and show that the linear span of functions obtained in this way
together with 
$\cA_{[n]}$ is the subspace  $\cA_{[n+1]}$.

We start with the  operators $\partial/\partial z$ and $z \partial/\partial z$,
\begin{align}
\frac{\partial}{\partial z} \frac{z^p\,\ov z^q}{(x-y)^n}= \frac{pz^{p-1}\ov z^q}{(x-y)^n}
-\frac{nz^p\,\ov z^q}{(x-y)^{n+1}};
\\
z\frac{\partial}{\partial z} \frac{z^p\,\ov z^q}{(x-y)^n}= \frac{pz^{p}\ov z^q}{(x-y)^n}
-\frac{nz^{p+1}\,\ov z^q}{(x-y)^{n+1}}.
\end{align}
In both lines, the first summand is contained in
 $\cA_n$. Therefore the second summand is contained in 
$\wt \cA$. Therefore, for $(p,q)\in\Omega_n$ we have 
\begin{equation}
\frac{z^{p}\,\ov z^q}{(x-y)^{n+1}},\quad  \frac{z^{p+1}\,\ov z^q}{(x-y)^{n+1}},\quad
\frac{z^{p}\,\ov z^{q+1}}{(x-y)^{n+1}}\in\wt \cA
\label{eq:inwtA}
\end{equation}
\begin{figure}
$${\mathrm a)}  \epsfbox{omega.2}\qquad {\mathrm b)}  \epsfbox{omega.3}\qquad
{\mathrm c)}  \epsfbox{omega.4}$$
\caption{To the proof of Lemma \ref{l:1}.
	\newline
	a) The sets $\Omega_n\subset \Omega_{n+1}$.
	\newline
	b) To formula (\ref{eq:inwtA}).
	\newline
	c) To formula (\ref{eq:}).
	\label{fig:2}
}
\end{figure}
\noindent
(the last inclusion follows from the application of the operator
 $\ov z\partial/\partial \ov z$).
 These fractions cover the whole domain
 $\Omega_{n+1}$ except two points 
$(n+1,n)$, $(n,n+1)$, see Fig.  \ref{fig:2}.b.

Next,
$$
z^2\frac{\partial}{\partial z} \frac{z^p\,\ov z^q}{(x-y)^n}= \frac{pz^{p+1}\ov z^q}{(x-y)^n}
-\frac{nz^{p+2}\,\ov z^q}{(x-y)^{n+1}}.
$$
Here we must examine 3 cases.

The first case. Let $p< n$, $(p,q)\ne (n-1,n)$. Then the first summand
is an element  $\cA_n$, therefore the second summand  is contained in
 $\wt\cA$.  However, this inclusion was obtained above 
(\ref{eq:inwtA}).

The second case.
 Let $p=n$. In the right-hand side we get 
\begin{equation}
\frac{-nx^{n+1} y^{q+1}}{(x-y)^{n+1}}.
\label{eq:}
\end{equation}
Therefore this fraction is contained in
 $\wt \cA$. For  $q<n-1$ we get a fraction from the list
 (\ref{eq:inwtA}), but for
  $q=n-1$ we land to the point  $(n+1,n)\in \Omega_{n+1}$,
see Fig. \ref{fig:2}.c. Applying the operator 
$\ov z\partial/\partial \ov z$ we land to $(n,n+1)\in \Omega_{n+1}$.
Thus
$$\cA_{n+1}\subset \wt A.$$

3) Let $(p,q)=(n-1,n)$. Thus in the right-hand side we get 
$$
-\frac{x^{n+1} y^n}{(x-y)^{n+1}}+\frac{(n-1)x^ny^{n+1}}{(x-y)^{n+1}}
$$
We landed to both summands above. 
On the other hand we again get an element of 
 $\cA_{n+1}$. \hfill $\square$

\sm 

{\sc Proof of Theorem \ref{th:2}.} Consider the algebra  $\wt \cB$,
generated by the operators (\ref{eq:generators}).
The commutator of $a(z)\partial/\partial z$ with an operator
of multiplication by a function $\phi(z,\ov z)$
is
\begin{equation}
\Bigl[ a(z)\frac{\partial}{\partial z},\phi(z,\ov z) \Bigr]= 
a(z)\frac{\partial \phi(z,\ov z)}{\partial z}.
\label{eq:[]}
\end{equation}
This observation and the lemma imply
that operators of multiplications by elements
of the algebra   $\cA$ are contained in
$\wt \cB$. Therefore $\wt\cB$ contain all products of the form 
   (\ref{eq:cB}), i.~e., $\wt \cB\supset \cB$. 

It remains to show that
 $\cB$ is an algebra. Consider a product of two elements
 of 
$\cB$, i.~e., of two products of the form  (\ref{eq:cB}).
We can transpose functional factors with operators (\ref{eq:6-operatorov})
using (\ref{eq:[]}), in this way we can move a functional factor to the beginning in all summands
 of the expression.
Next, the operators 
$\partial/\partial z$, $z \partial/\partial z$, $z^2\partial/\partial z$
form a Lie algebra 
 $\frs\frl(2)$ with respect to commutation, the same remark take place for operators
 $\partial/\partial \ov z$, $\ov z \partial/\partial \ov z$, 
$\ov z^2\partial/\partial\ov z$. Thus we get a direct sum of two
Lie algebras  
 $\frs\frl(2)$ (the Lie algebra of operators (\ref{eq:corr01})--(\ref{eq:corr03}))
is the diagonal in this direct sum). Transposing factors
we can put them to the order (\ref{eq:cB}).
\hfill $\square$

\noindent
 Math. Dept., University of Vienna; \\
 Institute for Theoretical and Experimental Physics (Moscow); \\
 MechMath Dept., Moscow State University;\\
 Institute for Information Transmission Problems;\\
 yurii.neretin@univie.ac.at;
 URL: http://mat.univie.ac.at/$\sim$neretin/

\noindent

\end{document}